\documentclass[12pt,english]{article}
\usepackage[T1]{fontenc}
\usepackage[latin9]{inputenc}
\usepackage{geometry}
\usepackage{amsmath}
\usepackage{amssymb}
\usepackage{stackrel}
\usepackage{setspace}
\makeatletter
\@ifundefined{date}{}{\date{}}
\makeatother

\usepackage{babel}
\begin{document}
\title{Some Remarks on the Notion of Bohr Chaos and Invariant Measures}
\author{Matan Tal\\
\\
The Hebrew University of Jerusalem}
\maketitle
\begin{abstract}
The notion of Bohr chaos was introduced in \cite{key-5,key-4}. We
answer a question raised in \cite{key-5} of whether a non uniquely
ergodic minimal system of positive topological entropy can be Bohr
chaotic. We also prove that all systems with the specification property
are Bohr chaotic, and by this line of thought give an independent
proof (and stengthening) of theorem 1 of \cite{key-5}. In addition,
we present an obstruction for Bohr chaos: a system with fewer than
a continuum of ergodic invariant probability measures cannot be Bohr
chaotic.\\
\end{abstract}

\section{Introduction and Main Results}

The notion of Bohr Chaos was introduced recently in \cite{key-5,key-4},
but let us recall the definitions. Two bounded complex sequences $\left(a_{n}\right)_{n\geq0},\,\left(b_{n}\right)_{n\geq0}$
are said to be $orthogonal$ (or $uncorrelated$ \footnote{Actually it corresponds to the notion of correlation only if one of
the sequences has such an average equal to $0$. But as we shall see
in a moment that is the only interesting case.}) if $\frac{1}{N}\stackrel[n=0]{N-1}{\sum}a_{n}b_{n}\xrightarrow[N\rightarrow\infty]{}0$.
Given a topological dynamical system $\left(X,T\right)$, i.e. a compact
metric space $X$ equipped with a continuous transformation $T:X\rightarrow X$,
a bounded real sequence $\left(a_{n}\right)_{n\geq0}$ is said to
be $orthogonal$ to it if it is orthogonal to all sequences $\left(f\left(T^{n}x\right)\right)_{n\geq0}$
where $x\in X$ and $f$ is a real valued continuous function. Last,
$\left(X,T\right)$ is said to be $Bohr\,chaotic$ if it is not orthogonal
to any non-trivial bounded real sequence $\left(a_{n}\right)_{n\geq0}$,
namely that there exists $x\in X$ and $f\in C\left(X\right)$ for
which $\left(f\left(T^{n}x\right)\right)_{n\geq0}$ correlates with
$\left(a_{n}\right)_{n\geq0}$ in the sense: $\limsup_{N\rightarrow\infty}\frac{1}{N}\left|\stackrel[n=0]{N-1}{\sum}a_{n}f\left(T^{n}x\right)\right|>0$.
By non-triviality of $\left(a_{n}\right)_{n\geq0}$ we mean $\limsup_{N\rightarrow\infty}\frac{1}{N}\stackrel[n=0]{N-1}{\sum}\left|a_{n}\right|>0$.\\

Notice that defining Bohr chaos with complex sequences and functions
will yield an equivalent definition since any complex sequence $\left(a_{n}\right)_{n\geq0}$
with $\limsup_{N\rightarrow\infty}\frac{1}{N}\stackrel[n=0]{N-1}{\sum}\left|a_{n}\right|>0$
satisfies 

$\limsup_{N\rightarrow\infty}\frac{1}{N}\stackrel[n=0]{N-1}{\sum}\left|\Re\left(a_{n}\right)\right|>0$
or $\limsup_{N\rightarrow\infty}\frac{1}{N}\stackrel[n=0]{N-1}{\sum}\left|\Im\left(a_{n}\right)\right|>0$.\\

One can define an a priori stronger notion by removing the boundedness
requirement on the sequence $\left(a_{n}\right)_{n\geq0}$. It seems
to be unknown whether the two notions are equivalent or not.\\

Throughout the paper, by a measure preserving system we will mean
a topological dynamical system (see the definition above) together
with an invariant probability measure.\\

Recall that two measure preserving systems $\left(X,T,\mu\right)$,
$\left(Y,S,\nu\right)$ are measure-theoretically disjoint, by definition,
if the only $T\times S$ invariant measure on $X\times Y$ that projects
onto $\mu$ and $\nu$ respectively is the product measure $\mu\times\nu$
(this definition originates in \cite{key-1-1}). In Section 2 we observe
that a Bohr chaotic topological dynamical system $\left(X,T\right)$
must have for any ergodic dynamical system $\left(Y,S,\nu\right)$
an invariant measure for which the two systems will not be disjoint,
and prove some stronger statements along this line. \\

Given a uniquely ergodic system, by taking $a_{n}=\xi^{n}$ where
$\xi\in S^{1}=\left\{ z\in\mathbb{C}\,:\,\left|z\right|=1\right\} $
is not one of the eigenvalues of its Koopman operator on $L^{2}$
(defined relative to its invariant probability measure), it is shown
in \cite{key-5} that it is not Bohr chaotic. Section 2 ends with
the following corollary: any system possessing a cardinality strictly
smaller than $2^{\aleph_{0}}$ of ergodic invariant probability measures
cannot be Bohr Chaotic (Cor. 2.5). \\

By the well-known measure-theoretic disjointness of a system of zero
measure-theoretic entropy to the Bernoulli-shift with two symbols
(since it is a system of completely positive entropy), it is shown
in \cite{key-5} that any topological dynamical system of zero topological
entropy is not Bohr chaotic. Thus a topological dynamical system cannot
be Bohr chaotic if it is either of zero topological entropy or uniquely
ergodic. Hence, the authors of \cite{key-5} implicitly raised the
question about the possibility of Bohr chaos in $minimal$ systems
(i.e. systems in which every orbit is dense) of positive topological
entropy which are not uniquely ergodic (by Section 2 here we know
they should actually have a continuum of ergodic invariant probability
measures). In Section $3$ we present an example of a minimal system
which is indeed Bohr chaotic. The example's construction is inspired
by the example of B. Weiss of a minimal system which is $universal$,
in the sense that for any aperiodic probability measure preserving
system it possesses an invariant measure which makes it measure-theoretically
isomorphic to that system (as usual, modulo null sets). So, in particular,
Weiss' system has a continuum of ergodic invariant measures. \\

In Section 4 we prove that invertible topological dynamical systems
admitting the specification property are Bohr chaotic. Such systems
are known to be universal for probability measure preserving systems
with entropy less than theirs (see \cite{key-8}). It is not known
to us whether such a universality is by itself a sufficient condition
for a topological dynamical system to be Bohr chaotic. In the end
of Section 4, we show how these results imply an independent proof
and a generalization of Theorem 1 of \cite{key-5}.\\

All sections are independent of one another although they share the
common thread of linking Bohr chaos with an abundance of invariant
measures.\\

\textbf{Acknowledgements.} The author wishes to thank N. Chandgotia,
A. Fan, H. Furstenberg, T. Meyerovitch and B. Weiss for their useful
advice. This research was partially supported by the Israeli Science
Foundation grant No. 1052/18.\\

\section{Bohr Chaos, Invariant Measures and Spectra}

In \cite{key-5} it is proved that a uniquely ergodic topological
dynamical system is not Bohr Chaotic. Arguing along similar lines,
it is possible to prove a more general fact. Here $C\left(Z\right)$
denotes complex valued functions defined on a compact metric space
$Z$.\\

\textbf{Lemma 2.1: }Let $\left(Y,S\right)$ be a topological dynamical
system with invariant measure $\mu$, $y_{0}$ a generic point for
$\mu$ and $g\in C\left(Y\right)$. If $\left(X,T\right)$ is a topological
dynamical system, and

$\limsup_{N\rightarrow\infty}\frac{1}{N}\left|\stackrel[n=0]{N-1}{\sum}f\left(T^{n}x_{0}\right)g\left(S^{n}y_{0}\right)\right|>0$
for some $x_{0}\in X$ and $f\in C\left(X\right)$, then there exits
a $T$-invariant measure $\eta$ and a joining of $\left(X,T,\eta\right)$
and $\left(Y,S,\mu\right)$ for which the integral of $f\left(x\right)g\left(y\right)$
does not vanish. If, in addition, $\mu$ is assumed to be ergodic,
then there exist such an $\eta$ and joining that are also ergodic.\\

\textbf{Proof:} There exists a monotone sequence $N_{k}$ such that
$\underset{k\rightarrow\infty}{\lim}\frac{1}{N_{k}}\left|\stackrel[n=0]{N_{k}-1}{\sum}f\left(T^{n}x_{0}\right)g\left(S^{n}y_{0}\right)\right|$
exists and is not zero. 

Passing to a subsequence if necessary $\frac{1}{N_{k}}\stackrel[n=0]{N_{k}-1}{\sum}\delta_{T^{n}x_{0}}\times\delta_{S^{n}y_{0}}$
converges to a measure $\theta$ on $\left(X\times Y,T\times S\right)$.
Taking $\eta$ to be the projection of $\theta$ onto $X$, $\theta$
becomes a joining $\lambda$ of $\left(X,T,\eta\right)$ and $\left(Y,S,\mu\right)$,
and $\underset{X\times Y}{\int}f\left(x\right)g\left(y\right)\,d\theta\left(x,y\right)=\underset{k\rightarrow\infty}{\lim}\frac{1}{N_{k}}\stackrel[n=0]{N_{k}-1}{\sum}f\left(T^{n}x_{0}\right)g\left(S^{n}y_{0}\right)\neq0$.\\

To prove the last statement of the lemma, we now assume $\mu$ to
be ergodic. The ergodic decomposition of $\lambda$ decomposes it
into ergodic joinings of ergodic $T$-invariant measures on $X$ and
the ergodic $S$-invariant measure $\mu$, and on a positive $\lambda$-measure
of them the integral of $f\left(x\right)g\left(y\right)$ does not
vanish.$\blacksquare$

\textbf{Theorem 2.2:} If $\left(X,T\right)$ is Bohr chaotic and $\left(Y,S,\mu\right)$
is an ergodic system then for every $g\in C\left(Y\right)$ that does
not vanish $\mu$-a.s., there exists a function $f\in C\left(X\right)$,
an ergodic $T$-invariant measure $\eta$ and an ergodic joining of
$\left(X,T,\eta\right)$ and $\left(Y,S,\mu\right)$ for which the
integral of $f\left(x\right)g\left(y\right)$ does not vanish. (So,
in particular, if $\left(Y,S,\mu\right)$ is non-trivial then the
two systems are not disjoint.)\\

\textbf{Proof:} Take $y_{0}\in Y$ which is a generic point of $\left(Y,S,\mu\right)$,
then $\lim_{N\rightarrow\infty}\frac{1}{N}\stackrel[n=0]{N}{\sum}\left|g\left(S^{n}y_{0}\right)\right|>0$.
Bohr chaoticity then implies that there exists $x_{0}\in X$ and $f\in C\left(X\right)$
for which

$\limsup_{N\rightarrow\infty}\frac{1}{N}\left|\stackrel[n=0]{N-1}{\sum}f\left(T^{n}x_{0}\right)g\left(S^{n}y_{0}\right)\right|>0$
as needed in order to apply Lemma 2.1. $\blacksquare$\\

It is unknown whether the converse of Theorem 2.2 is true: If for
every $\left(Y,S,\mu\right)$ and $g$ as above there exist $f,\eta$
and joining as above, does this imply that $\left(X,T\right)$ is
Bohr chaotic? Even the answer to the following natural question is
not known to be positive: in order to detect Bohr chaoticity of $\left(X,T\right)$,
is it sufficient to consider sequences $\left(a_{n}\right)_{n\geq0}\in\mathbb{\overline{D}}^{\mathbb{N}\cup\left\{ 0\right\} }$
which are generic points of the shift space $\mathbb{\overline{D}}^{\mathbb{N}\cup\left\{ 0\right\} }$
for ergodic shift-invariant measures, where $\mathbb{\overline{D}}$
is the closed unit disk in the complex plane? A postive answer to
the first question implies a positive answer to the second.\\

\textbf{Lemma 2.3:} For every ergodic joining $\theta$ of a measure
preserving system $\left(X,T,\mu\right)$ and a circle rotation by
$\xi\in S^{1}=\left\{ z\in\mathbb{C}\,:\,\left|z\right|=1\right\} $
such that $\xi$ is not an eigenvalue of the first system's Koopman
operator, and every $f\in L^{1}\left(X,\mu\right)$, the integral
of $f\left(x\right)e^{2\pi iy}$ with respect to $\theta$ vanishes.\\

\textbf{Proof:} Assume the integral does not vanish, then, by the
Mean Ergodic Theorem, there exists a monotone sequence $N_{k}$ for
which $e^{2\pi y}\cdot\frac{1}{N_{k}}\stackrel[n=0]{N_{k}-1}{\sum}f\left(T^{n_{k}}x\right)\xi^{n}$
converges to the value of the integral $\theta$-a.s.. Because this
value is not zero this implies that the conditional measure of $\theta$
relative to the $X$ fibers is $\mu$-a.s. supported on just one point
of $Y$ which we shall denote $g\left(x\right)$. The invariance of
$\theta$ implies that $g\left(Tx\right)=\xi g\left(x\right)$ $\mu$-a.s.,
contradicting the fact that $\xi$ is not an eigenvalue of $\left(X,T,\mu\right)$.
$\blacksquare$\\

\textbf{Theorem 2.4:} If a topological dynamical system $\left(X,T\right)$
is Bohr Chaotic then for every $\xi\in S^{1}$ there exists an ergodic
$T$-invariant measure such that $\xi$ is an eigenvalue of the Koopman
operator on $L^{2}\left(X,\mu\right)$.\\

\textbf{Proof:} For $R_{\xi}:S^{1}\rightarrow S^{1}$ the rotation
by $\xi$ and $m$ the Haar measure on $S^{1}$, we can apply Theorem
2.2 with $\left(Y,S,\mu\right)=\left(S^{1},R_{\xi},m\right)$ and
$g\left(y\right)=e^{2\pi iy}$ and obtain a function $f\in C\left(X\right)$,
an ergodic $T$-invariant measure $\eta$ and an ergodic joining of
$\left(X,T,\eta\right)$ and $\left(S^{1},R_{\xi},m\right)$ such
that the integral of $f\left(x\right)e^{2\pi iy}$ does not vanish.
By Lemma 2.3, $\xi$ is an eigenvalue as required. $\blacksquare$\\

The converse to Theorem 2.4 is false. As a counter-example consider
the sub-shift of $\left\{ 0,1\right\} ^{\mathbb{Z}}$ generated by
the union of all Sturmian sub-shifts (it is the closure of this union).
It satisfies the consequent of Theorem 2.4 regarding the spectra of
the invariant probability measures, but it is not Bohr chaotic because
it has zero topological entropy (in any point of this system appear
no more than $N+1$ different words of length $N$, and so it is immediate
that any ergodic measure has zero entropy by considering one of its
generic points).\\

\textbf{Corollary 2.5:} If a topological dynamical system $\left(X,T\right)$
possesses a cardinality strictly smaller than $2^{\aleph_{0}}$ of
ergodic invariant measures then it is not Bohr Chaotic.\\

\textbf{Proof :} If $\left\{ \mu_{i}\right\} _{i}$ are the ergodic
invariant measures of the system then

$\Sigma\left(\mu_{i}\right):=\left\{ \xi\in S^{1}\,:\,\xi\,is\,an\,eigenvalue\,of\,the\,Koopman\,operator\,on\,L^{2}\left(X,\mu_{i}\right)\right\} $
is countable and hence the cardinality of $\Sigma=\underset{i}{\cup}\Sigma\left(\mu_{i}\right)$
is strictly smaller than $2^{\aleph_{0}}$, so there must be some
$\xi\in S^{1}$ for which $\xi\in S^{1}\setminus\Sigma$ . By Theorem
2.4 this implies the desired result. $\blacksquare$\\

\section{An Example of a Minimal System Which is Bohr Chaotic}

The idea for the following construction is taken from B. Weiss' example
of a universal minimal system given in \cite{key-9} .\\

\textbf{Lemma 3.1:} Assume $L>2$, $X\subseteq\text{\ensuremath{\Lambda}}^{\mathbb{Z}}$
($\Lambda$ a finite set) is a mixing shift of finite type (hereinafter
abbreviated as SFT), $W$ an allowable word in it and $K>0$ such
that for any allowable words $\alpha,\beta$ and $l\geq K$ there
exists an allowable word $\alpha w\beta$ such that $w$ is a word
of length $l$ that contains $W$ (such a $K$ always exists in a
mixing SFT). Then the non-empty SFT $X_{L,W,K}\subseteq X$ composed
of the points in $X$ for which every word of length $LK$ contains
$W$ is mixing.\\

\textbf{Proof:} We construct for any $\left(x_{n}\right)_{n\in\mathbb{Z}},\,\left(y_{n}\right)_{n\in\mathbb{Z}}\in X_{L,W,K}$
a $\left(z_{n}\right)_{n\in\mathbb{Z}}\in X_{L,W,K}$ such that $z_{n}=x_{n}$
for $n<0$ and $z_{n}=y_{n}$ for $n\geq LK$. First, of course we
define $z_{n}=x_{n}$ for $n<0$ and $z_{n}=y_{n}$ for $n\geq LK$.
If $W$ appears in $x_{-length\left(W\right)+1}\dots x_{0}\dots x_{K-1}$
then we define $z_{0}\dots z_{K-1}:=x_{0}\dots x_{K-1}$, otherwise
we take $z_{0}\dots z_{K-1}$ to be any word that contains $W$ such
that $\dots x_{-3}x_{-2}x_{-1}z_{0}\dots z_{K-1}$ is allowable in
$X$. And if $W$ appears in $y_{\left(L-1\right)K}\dots y_{LK-1}\dots y_{LK+length\left(W\right)-2}$
then we define

$z_{\left(L-1\right)K}\dots z_{LK-1}:=y_{\left(L-1\right)K}\dots y_{LK-1}$,
otherwise we take $z_{\left(L-1\right)K}\dots z_{KL-1}$ to be any
word that contains $W$ such that 

$z_{\left(L-1\right)K}\dots z_{KL-1}y_{KL}y_{KL+1}y_{KL+2}\dots$
is allowable in $X$. $L>2$ and so it remains to take any $z_{K}\dots z_{\left(L-1\right)K-1}$
that contains $W$ and for which $\left(z_{n}\right)_{n\in\mathbb{Z}}\in X$.
It follows that $\left(z_{n}\right)_{n\in\mathbb{Z}}\in X_{L,W,K}$.
$\blacksquare$\\

\textbf{Lemma 3.2:} Assume we are given a sequence of positive integers
$L_{i}$ satisfying $\stackrel[i=1]{\infty}{\sum}\frac{16}{L_{i}}<1$.
With the aid of Lemma 3.1, we define inductively a decreasing sequence
of mixing SFTs:

$X_{0}:=\left\{ -1,1\right\} ^{\mathbb{Z}}$ and $X_{i}$ is the mixing
SFT with $K_{i}$ in the role of $K$ in that lemma when applied to
$X_{i-1}$, $L_{i}$ and $W_{i-1}$ - some allowable word of $X_{i-1}$
. Then the sub-shift $X=\stackrel[i=0]{\infty}{\cap}X_{i}\neq\emptyset$
is Bohr chaotic.\\

\textbf{Proof:} Let $\left(a_{n}\right)_{n\in\mathbb{Z}}$ be a bounded
real sequence satisfying $\lim_{k\rightarrow\infty}\frac{1}{N_{k}}\stackrel[n=0]{N_{k}-1}{\sum}\left|a_{n}\right|>0$.
And let us assume without loss of generality that $N_{k}>2N_{k-1}$
for all $k$. We define a sequence $\left(x_{n}^{\left(0\right)}\right)_{n\in\mathbb{Z}}$
by
\[
x_{n}^{\left(0\right)}=\left\{ \begin{array}{cc}
1 & a_{n}\geq0\\
-1 & a_{n}<0
\end{array}\right..
\]
\\
Then $\stackrel[n=0]{N_{k}-1}{\sum}a_{n}x_{n}^{\left(0\right)}=\stackrel[n=0]{N_{k}-1}{\sum}\left|a_{n}\right|$
for all $k$. The point is we can recursively modify $\left(x_{n}^{\left(i-1\right)}\right)_{n\in\mathbb{Z}}$
to form $\left(x_{n}^{\left(i\right)}\right)_{n\in\mathbb{Z}}\in X_{i}$
(for $i=1,2,\dots$) such that for all $k$
\[
\left(*\right)\,\,\stackrel[n=0]{N_{k}-1}{\sum}a_{n}x_{n}^{\left(i\right)}>\left(1-\stackrel[j=1]{i}{\sum}\frac{16}{L_{j}}\right)\cdot\stackrel[n=0]{N_{k}-1}{\sum}\left|a_{n}\right|.
\]
 This implies that any of $\left(x_{n}^{\left(i\right)}\right)_{n\in\mathbb{Z}}$'s
partial limits (as a sequence in $i$) $\left(x_{n}\right)_{n\in\mathbb{Z}}$
satisfies $\stackrel[n=0]{N_{k}-1}{\sum}a_{n}x_{n}>\left(1-\stackrel[j=1]{\infty}{\sum}\frac{16}{L_{j}}\right)\cdot\stackrel[n=0]{N_{k}-1}{\sum}\left|a_{n}\right|$
for all $k$ and thus 
\[
\liminf_{k\rightarrow\infty}\frac{1}{N_{k}}\stackrel[n=0]{N_{k}-1}{\sum}a_{n}x_{n}\geq\left(1-\stackrel[j=1]{\infty}{\sum}\frac{16}{L_{j}}\right)\cdot\lim_{k\rightarrow\infty}\frac{1}{N_{k}}\stackrel[n=0]{N_{k}-1}{\sum}\left|a_{n}\right|>0.
\]
It remains to show how those $\left(x_{n}^{\left(i\right)}\right)_{n\in\mathbb{Z}}\in X_{i}$
fulfilling condition $\left(*\right)$ are constructed. For simplicity
of notation let us assume that every $L_{i}$ is a multiple of $8$.\\

For every integer $l$, we divide the segment $x_{lK_{i}\frac{L_{i}}{2}}^{\left(i-1\right)},\dots,x_{\left(l+1\right)K_{i}\frac{L_{i}}{2}-1}^{\left(i-1\right)}$
into disjoint consecutive words of length $K_{i}$: $V_{1},\dots,V_{\frac{L_{i}}{2}}$.
If, for every $l$, we modify one of these $V_{i}$ (say, sequentialy
for the non-negative $l$s and then sequentialy for the negative ones)
such that it will contain $W_{i}$ while staying in $X_{i-1}$, then
the new point will lie in $X_{i}$. We now describe how to choose
a $V_{i}$ to modify for every $l$ such that also condition $\left(*\right)$
will be satisfied. \\

If there is no $k$ for which $lK_{i}\frac{L_{i}}{2}\leq N_{k}<\left(l+1\right)K_{i}\frac{L_{i}}{2}$
then $\left(x_{n}^{\left(i\right)}\right)_{n\in\mathbb{Z}}$ is defined
by selecting from each such sequence of words $V_{1},\dots,V_{\frac{_{L_{i}}}{2}}$
one for which the sum of $\left|a_{n}\right|$ on its indices is minimal
and modifying it to contain $W_{i}$ - this is possible by the choice
of $K_{i}$.\\

Else, if $l>0$ there is exactly one $k_{0}$ for which $lK_{i}\frac{L_{i}}{2}\leq N_{k_{0}}<\left(l+1\right)K_{i}\frac{L_{i}}{2}$
(since $N_{k}>2N_{k-1}$ for all $k$), so this index $N_{k_{0}}$
belongs either to $V_{1},\dots,V_{\frac{_{L_{i}}}{4}}$ or to $V_{\frac{L_{i}}{4}+1},\dots,V_{\frac{_{L_{i}}}{2}}$.
In the one of these two word collections of which $N_{k_{0}}$ does
not belong, we select the word for which the sum of $\left|a_{n}\right|$
on its indices is minimal and modify it to contain $W_{i}$.\\

The case which is left to handle is the case where $l=0$ and there
exist $k$s for which $0<N_{k}<K_{i}\frac{L_{i}}{2}$. But then again,
although there can be many $k$s for which $0<N_{k}<K_{i}\frac{L_{i}}{4}$,
there can be at most one $k_{0}$ for which $K_{i}\frac{L_{i}}{4}\leq N_{k_{0}}<K_{i}\frac{L_{i}}{2}$.
So such an index $N_{k_{0}}$ can belong either to at most one of
$V_{\frac{L_{i}}{4}+1},\dots,V_{\frac{_{3L_{i}}}{8}}$ or to $V_{\frac{_{3L_{i}}}{8}+1},\dots,V_{\frac{_{L_{i}}}{2}}$.
And again, in the one of these two word collections of which $N_{k_{0}}$
does not belong, we select the word for which the sum of $\left|a_{n}\right|$
on its indices is minimal and modify it to contain $W_{i}$. $\blacksquare$
\\

\textbf{Theorem 3.3:} There exists a minimal topological dynamical
system that is Bohr Chaotic.\\

\textbf{Proof:} We will construct a minimal sub-shift $X\subseteq\left\{ -1,1\right\} ^{^{\mathbb{Z}}}$
that satisfies the conditions of Lemma 3.2. First let us fix some
sequence of positive integers $L_{i}$ satisfying $\stackrel[i=1]{\infty}{\sum}\frac{16}{L_{i}}<1$.\\

We construct a sub-shift $X$ of the shift system $\left(\left\{ -1,1\right\} ^{^{\mathbb{Z}}},\sigma\right)$
as the intersection of a descending family of mixing SFTs $X_{0}\supseteq X_{1}\supseteq X_{2}\supseteq\dots$
defined recursively (the recursion will also define $0<K_{1}<K_{2}<\dots$
and $W_{0},W_{1},\dots\in\underset{n\in\mathbb{N}}{\cup}\left\{ -1,1\right\} ^{n}$
which are certain allowable words in the corresponding $X_{i}$s):
\begin{itemize}
\item $X_{0}$ will be the full shift: all possible sequences of $-1$s
and $1$s. It is of course a mixing SFT, $"-1\,1"$ is an allowable
word and we denote it by $W_{0}$.\\
\item Since $X_{i-1}$ is a mixing SFT there exists $K_{i}>0$ such that
for any allowable of words $\alpha,\beta$ and $l\geq K_{i}$ there
exists an allowable word $\alpha w\beta$ such that $w$ is a word
of length $l$ that contains $W_{i-1}$. $X_{i}$ will be the sub-shift
of $X_{i-1}$ for which every word of length $L_{i}K_{i}$ contains
$W_{i-1}$. By lemma 3.1 it is a mixing SFT.

Choose $W_{i}$ to be a word containing all $X_{i}$ allowable words
of length $i+1$.
\end{itemize}
$X:=\underset{i}{\cap}X_{i}$ is non-empty by compactness of each
$X_{i}$, and it is minimal because by the choice of the $W_{i}$s
all points have dense orbits. $\blacksquare$\\

\section{The Specification Property implies Bohr Chaos}

\subsection{Proof for a Sub-Shift}

Let $\Lambda$ be a finite set, $X$ a sub-shift of the shift system
$\left(\Lambda^{\mathbb{N}\cup\left\{ 0\right\} },\sigma\right)$
containing more than one point in its eventual image (i.e. $\stackrel[i=1]{\infty}{\cap}\sigma^{i}\left(X\right)$).
We shall say it satisfies $the\,specification\,property$, if there
exists some $u\geq0$ such that for any two allowable words in the
eventual image (of any length) $\alpha$ and $\beta$ there exists
a point $x\in\stackrel[i=1]{\infty}{\cap}\sigma^{i}\left(X\right)$
with $\alpha$ and $\beta$ appearing in it spaced by exactly $u$
coordinates (on which there is no constraint) \footnote{This definition of the specification property for systems with an
eventual image containing more than one point is weaker than the usual
one for such systems (cf. \cite{key-1-2} for the usual one).}. Mixing sofic sub-shifts containing more than one point are such
examples. Notice that all sub-shifts admitting specification have
positive topological entropy. In the present section we prove that
such an $X$ is Bohr chaotic. \\

\textbf{Theorem 4.1.1:} A sub-shift $X$ admitting specification is
Bohr chaotic.\\

\textbf{Proof}: It is enough to prove that $\left(\stackrel[i=1]{\infty}{\cap}\sigma^{i}\left(X\right),\sigma\right)$
is Bohr chaotic. The idea is to take two different allowable words
in the eventual image $\alpha$ and \textbf{$\beta$} of equal length
$l$, and to consider the sub-shift $Y\subseteq\stackrel[i=1]{\infty}{\cap}\sigma^{i}\left(X\right)$
containing exactly the points made of concatenations of infinitely
many blocks which are either $\alpha$ or $\beta$ and with spaces
$u$ between consecutive blocks (the spaces can contain anything).
Now, if all such points could be locally parsed uniquely, then it
would be possible to define a $\sigma$-equivariant continuous map
$\varphi:Y\rightarrow\left(\Lambda\cup\left\{ 0\right\} \right)^{\mathbb{N}\cup\left\{ 0\right\} }$
that receives such a point and ``forgets'' the information in the
spaces: fills them all with zeros. Denoting by $Z$ the image of $\varphi$,
and by $Z\,'$ all points of $Z$ of which an $\alpha$ or a $\beta$
begins at the zero coordinate. $Z\,'$ is not invariant under $\sigma$,
however it is invariant under $\sigma^{l+u}$ and $\left(Z\,',\sigma^{l+u}\right)$
is naturally isomorphic to $\left(\left\{ 0,1\right\} ^{\mathbb{N}\cup\left\{ 0\right\} },\sigma\right)$
by an isomorphism we shall denote $\psi$ (say an $\alpha$ corresponds
to a $0$ and a $\beta$ to a $1$). $Z\,',\sigma\left(Z\,'\right),\dots,\sigma^{l+u-1}\left(Z\,'\right)$
are pairwise disjoint, and this easily implies that $\left(Z,\sigma\right)$
is Bohr chaotic. Now, to deduce that $\left(\stackrel[i=1]{\infty}{\cap}\sigma^{i}\left(X\right),\sigma\right)$
is Bohr chaotic, one can use Theorem 1 in \cite{key-5}, but there
is no need for it in this simple setting. For if $\limsup_{N\rightarrow\infty}\frac{1}{N}\stackrel[n=0]{N-1}{\sum}\left|a_{n}\right|>0$
then for some $0\leq t<l+u$ also $\limsup_{N\rightarrow\infty}\frac{1}{N}\stackrel[n=0]{N-1}{\sum}\left|a_{\left(l+u\right)n+t}\right|>0$.
So if $z=\left(z_{n}\right)_{n\in\mathbb{N}\cup\left\{ 0\right\} }$
is the element in $Z$ satisfying $\sigma^{t}z=\psi^{-1}\left(\left(sgn\left(a_{\left(l+u\right)n+t}\right)\right)_{n=0}^{\infty}\right)$,
and defining $f:Z\rightarrow\left\{ -1,0,1\right\} $ to be $0$ outside
of $Z'$ and equal $\pi_{0}\circ\psi$ on $Z'$ ($\pi_{0}$ being
the projection on the zero co-ordinate), then sequence $f\left(\sigma^{n}z\right)$
satisfies that $\limsup_{N\rightarrow\infty}\frac{1}{N}\stackrel[n=0]{N-1}{\sum}a_{n}f\left(\sigma^{n}z\right)>0$.\\

It is left to find words $\alpha$ and $\beta$ for which the parsing
will indeed be unique. Taking $\alpha$ and $\beta$ to be two words
as in Lemma 4.1.2 (applied for the sub-shift $\stackrel[i=1]{\infty}{\cap}\sigma^{i}\left(X\right)$)
for $r=\frac{1}{3}$, non-unique parsing implies $l-\frac{l}{3}-u<\frac{l}{3}$
which is equivalent to $l<3u$. So taking $l\geq3u$ guarantees that
the parsing will be unique. $\blacksquare$\\

We define the overlap between two different words $a_{0}\dots a_{r-1}$
and $b_{0}\dots b_{s-1}$ to be 
\[
\max\left(\left\{ 1\leq k\leq r,s\,:\,a_{r-k}\dots a_{r-1}=b_{0}\dots b_{k-1}\,or\,b_{s-k}\dots b_{s-1}=a_{0}\dots a_{k-1}\right\} \cup\left\{ 0\right\} \right),
\]
 and the self-overlap of a word $a_{1}\dots a_{m}$ to be 
\[
\max\left(\left\{ 1\leq k<m\,:\,a_{r-k}\dots a_{r-1}=a_{0}\dots a_{k-1}\right\} \cup\left\{ 0\right\} \right).
\]
\\

\textbf{Lemma 4.1.2:} Let $X$ be a sub-shift of topological entropy
$h>0$. For every $0<r<1$ and any large enough $l$ there exists
two different allowable words $\alpha$ and $\beta$ of length $l$
such that the overlap between $\alpha$ and $\beta$ and the self-overlaps
of $\alpha$ and $\beta$ are all smaller than $rl$.\\

\textbf{Proof:} Given an ergodic invariant measure $\mu$ of entropy
$\tilde{h}$ on $\left(X,\sigma\right)$, by Theorem 1 of \cite{key-2},

$x=\left(x_{n}\right)_{n\in\mathbb{N}\cup\left\{ 0\right\} }\in X$
$\mu$-almost-surely satisfies $\underset{n\rightarrow\infty}{\lim}\frac{\log R_{n}\left(x\right)}{n}=\tilde{h}$
where 

$R_{n}\left(x\right)=\min\left\{ k>0\,:\,x_{k}\dots x_{k+n-1}=x_{0}\dots x_{n-1}\right\} $
(as one would probably guess from the Shannon-Macmillan-Breiman Theorem).
In fact, Theorem 1 of \cite{key-2} says exactly this but for a slightly
different definiton of $R_{n}\left(x\right)$. It is defined as $\min\left\{ k>n-1\,:\,x_{k}\dots x_{x_{k+n-1}}=x_{0}\dots x_{n-1}\right\} $,
and let us denote this definition as $R'_{n}\left(x\right)$. We claim
that for any given $x$ - unless $\underset{n\rightarrow\infty}{\liminf}\frac{\log R'_{n}\left(x\right)}{n}=0$
(and then also $\underset{n\rightarrow\infty}{\liminf}\frac{\log R_{n}\left(x\right)}{n}=0$)
- $R_{n}\left(x\right)$ and $R'_{n}\left(x\right)$ differ only for
a finite set of $n$-s (and hence

$\underset{n\rightarrow\infty}{\lim}\frac{\log R_{n}\left(x\right)}{n}=\underset{n\rightarrow\infty}{\lim}\frac{\log R'_{n}\left(x\right)}{n})$.
For if there exists an infinite strictly monotone increasing sequence
$\left(n_{k}\right)_{k=1}^{\infty}$ for which $R_{n_{k}}\left(x\right)\neq R'_{n_{k}}\left(x\right)$,
then $R'_{R_{n_{k}}\left(x\right)}\left(x\right)=R_{n_{k}}\left(x\right)$,
and hence

$\lim_{k\rightarrow\infty}\frac{\log R'_{R_{n_{k}}\left(x\right)}\left(x\right)}{R_{n_{k}}\left(x\right)}=0$.
So Theorem 1 of \cite{key-2} indeed implies the result stated at
the begining of this paragraph.\\

Taking $\mu$ to be an ergodic invariant measure of positive entropy
on $\left(X,\sigma\right)$, we obtain

$\underset{n\rightarrow\infty}{\lim}\frac{\log R_{n}\left(x\right)}{n}=h>0$
$\mu$-almost-surely. By Egorov's Theorem this implies that for any
large enough $l$ there is a positive $\mu$-probability that $x$
will satisfy $R_{\left\lfloor rl\right\rfloor }\left(x\right),R_{\left\lfloor rl\right\rfloor }\left(\sigma^{l-\left\lfloor rl\right\rfloor }x\right),R_{\left\lfloor rl\right\rfloor }\left(\sigma^{l}x\right)>2l$.
For any such $x$, $\alpha=x_{0}\dots x_{l-1}$ and $\beta=x_{l}\dots x_{2l-1}$
will do. $\blacksquare$ \\

\textbf{Remark:} An alternative path for the result of this section
instead of using Lemma 4.1.2 would be to use the result of A. Bertrand
that for symbolic systems the specification property implies the synchronization
property (cf. \cite{key-1}) and that the synchronization property
for such systems implies that they are coded systems (cf. \cite{key-3}).
The reader can find the definitions of these terms in \cite{key-1}.

\subsection{Proof for a General Invertible Topological Dynamical System}

Let us define the specification property for a general topological
dynamical system - not necessarily a sub-shift. Let $\left(X,d\right)$
be a compact metric space and and $\left(X,T\right)$ a topological
dynamical system with an evetual image (i.e. $\stackrel[i=1]{\infty}{\cap}T^{i}\left(X\right)$)
containing more than one point (again, our definition will be weaker
for such systems than the usual one found in \cite{key-1-2}).

For any $z\in X$, an orbit segment $z,Tz,\dots,T^{r-1}z$ is said
to $\varepsilon-shadow$ a sequence $z^{\left(0\right)},\dots,z^{\left(r-1\right)}\in X$
if 
\[
d\left(z,z^{\left(0\right)}\right),d\left(Tz,z^{\left(1\right)}\right),\dots,d\left(T^{r-1}z,z^{\left(r-1\right)}\right)<\varepsilon.
\]
$\left(X,T\right)$ is said to satisfy $the\,specification\,property$,
if for any $\varepsilon>0$ there exists some $u_{\varepsilon}\geq0$
such that for any finite number of orbit segments in the eventual
image 

\[
\left(x^{\left(1\right)},Tx^{\left(1\right)},\dots,T^{r_{1}-1}x^{\left(1\right)}\right),\dots,\left(x^{\left(n\right)},Tx^{\left(n\right)},\dots,T^{r_{n}-1}x^{\left(n\right)}\right)
\]
there exists a point $z\in\stackrel[i=1]{\infty}{\cap}T^{i}\left(X\right)$
for which 

$z,Tz,\dots,T^{r-1}z$ $\varepsilon$-shadows $x^{\left(1\right)},Tx^{\left(1\right)},\dots,T^{r_{1}-1}x^{\left(1\right)}$,

$T^{r_{1}+u_{\varepsilon}}z,Tz,\dots,T^{r_{1}+r_{2}+u_{\varepsilon}-1}z$
$\varepsilon$-shadows $x^{\left(2\right)},Tx^{\left(2\right)},\dots,T^{r_{2}-1}x^{\left(2\right)}$,

$T^{r_{1}+r_{2}+2u_{\varepsilon}}z,Tz,\dots,T^{r_{1}+r_{2}+r_{3}+2u_{\varepsilon}-1}z$
$\varepsilon$-shadows $x^{\left(3\right)},Tx^{\left(3\right)},\dots,T^{r_{3}-1}x^{\left(3\right)}$
etc.\\

We define the $\varepsilon$-overlap between two different orbit segments
$x,Tx,\dots,T^{r-1}x$ and $y,Ty,\dots,T^{s-1}y$ of points $x,y\in X$
and the $\varepsilon$-self-overlap an orbit segment $x,f\left(x\right),\dots,f^{r-1}\left(x\right)$
in the same manner as the parallel notions of Section 4.1 only with
the equalities inside the curly brackets replaced by requiring $\varepsilon$-shadowing
(for example, we do not require $T^{r-k}x\dots T^{r-1}x=y\dots T^{k-1}y$
but that $T^{r-k}x\dots T^{r-1}x$ $\varepsilon$-shadows $y\dots T^{k-1}y$).\\

\textbf{Theorem 4.2.1:} Given a compact metric space $\left(X,d\right)$,
a topological dynamical system $\left(X,T\right)$ admitting specification
is Bohr chaotic.\\

\textbf{Proof:} Considering an $\varepsilon$ as in Lemma 4.2.2 with
$r=\frac{1}{3}$ (the lemma is applied on $\stackrel[i=1]{\infty}{\cap}T^{i}\left(X\right)$),
we then take two orbit segments in the eventual image $x,Tx,\dots,T^{l-1}x$
and $y,Ty,\dots,T^{l-1}y$ with small $\varepsilon$-overlaps as guaranteed
to exist by it and also require $l\geq3u_{\frac{\varepsilon}{3}}$.\\

We denote by $Y\subseteq\stackrel[i=1]{\infty}{\cap}T^{i}\left(X\right)$
the sub-system containing all points $z\in\stackrel[i=1]{\infty}{\cap}T^{i}\left(X\right)$
such that $\dots,T^{-2}z,T^{-1}z,z,Tz,T^{2}z,\dots$ is composed of
$\frac{\varepsilon}{3}$-shadowings of the above two orbit segments
with $u_{\frac{\varepsilon}{3}}$ spaces between them which can contain
anything. These points can be uniquely parsed, and from here we can
argue as in the first paragraph of Theorem 4.1.1 to finish the proof.
$\blacksquare$ \\

\textbf{Remark:} Notice that all we need for the proof of Theorem
4.2.1 is that for some $\varepsilon$ and $l$ there will be two different
orbit segments in $X$ of length $l$ with an $\varepsilon$-overlap
and $\varepsilon$-self-overlaps smaller than $\frac{l}{3}$ which
``admit specification'' for $\text{\ensuremath{\frac{\varepsilon}{3}}}$
with spaces $u_{\frac{\varepsilon}{3}}\leq\frac{l}{3}$. By '``admitting
specification'' we mean that we can combine them as we wish (up to
$\frac{\varepsilon}{3}$-shadowings) with spaces of size $u_{\frac{\varepsilon}{3}}$
to form a complete orbit. So we do not need the entire specification
property in all its strength.\\

\textbf{Lemma 4.2.2:} Given $0<r<1$, a compact metric space $\left(X,d\right)$
and a topological dynamical system $\left(X,T\right)$ of topological
entropy $h>0$, then for any sufficiently small $\varepsilon>0$,
any large enough $l$ satisfies that there exist two different orbit
segments $x,Tx,\dots,T^{l-1}x$ and $y,Ty,\dots,T^{l-1}y$ with the
$\varepsilon$-overlap between them and their $\varepsilon$-self-overlaps
all smaller than $rl$. \footnote{A similar result can be found in Lemma 5.6 of \cite{key-8}. }\\

\textbf{Proof:} First we state a theorem from \cite{key-12}. Denote
$B_{n,\varepsilon}\left(x\right):=\left\{ y\in X\,:\,\forall0\leq i<n\,d\left(T^{i}y,T^{i}x\right)<\varepsilon\right\} $
(the $n$-th Bowen ball of radius $\varepsilon$ around $x$), $R_{n,\varepsilon}\left(x\right):=\min\left\{ k>0\,:\,T^{k}x\in B_{n,\varepsilon}\left(x\right)\right\} $.
If $\mu$ is an ergodic invariant measure of entropy $\tilde{h}$
on $\left(X,T\right)$, then $\underset{n\rightarrow\infty}{\liminf\frac{\log R_{n,\varepsilon}\left(x\right)}{n}}$
equals a constant $C_{\varepsilon}$ $\mu$-almost-surely. Theorem
3 of \cite{key-12} (together with remark (D) above it) says that
$\underset{\varepsilon\rightarrow0}{\lim}C_{\varepsilon}=\tilde{h}$.

Using this theorem in the same spirit as Theorem 1 of \cite{key-2}
in the last paragraph of the proof of Lemma 4.1.2 yields the required
result. $\blacksquare$\\

For a topological dynamical system $\left(X,T\right)$, Theorem $1$
in \cite{key-5} says that if it contains a closed subset $Y$ for
which there exists a $k>0$ such that $Y$ is $T^{k}$-invariant and
$\left(Y,T^{k}\right)\cong\left(\Lambda^{\mathbb{N}},\sigma\right)$
(for some finite set $\Lambda$) then it is Bohr chaotic. It is interesting
to notice that the proof of Theorem 4.2.1 - together with the remark
following it - can serve us for an independent ergodic-theoretic proof
of this fact in an even more general formulation.\\

\textbf{Theorem 4.2.3:} If a topological dynamical system $\left(X,T\right)$
contains a closed subset $Y$ for which there exists a $k>0$ such
that $Y$ is $T^{k}$-invariant and $\left(Y,T^{k}\right)$ admits
specification then $\left(X,T\right)$ is Bohr chaotic.\\

\textbf{Proof:} We show that the system $\left(\stackrel[i=0]{k-1}{\cup}T^{i}Y,T\right)$
is a topological dynamical system that fulfills the requirement of
the remark following the proof of Theorem 4.2.1, thus it is Bohr chaotic
and thus so is $\left(X,T\right)$.\\

The eventual image of $\left(\stackrel[i=0]{k-1}{\cup}T^{i}Y,T\right)$
has positive topological entropy (because the eventual image of $\left(Y,T^{k}\right)$
has), and so it fits Lemma 4.2.2 with some $\varepsilon_{0}>0$.

We take a $\delta>0$ for which $d\left(x,y\right),d\left(Tx,Ty\right),d\left(T^{2}x,T^{2}y\right),\dots,d\left(T^{k-1}x,T^{k-1}y\right)<\frac{\varepsilon_{0}}{3}$
for any

$x,y\in\stackrel[i=0]{k-1}{\cup}T^{i}Y$ with $d\left(x,y\right)<\delta$.
\\

By Lemma 4.2.2 applied to the eventual image of $\left(\stackrel[i=0]{k-1}{\cup}T^{i}Y,T\right)$,
there exist in it two different orbit segments $x^{\left(1\right)},Tx^{\left(1\right)},\dots,T^{\left(l'-1\right)k}x^{\left(1\right)}$
and

$x^{\left(2\right)},Tx^{\left(2\right)},\dots,T^{\left(l'-1\right)k}x^{\left(2\right)}$
that have an $\varepsilon_{0}$-overlap and $\varepsilon_{0}$-self-overlaps
all smaller than $\frac{l'}{4}$ and $l'\geq3u_{\delta}$, where $u_{\delta}$
is the the specification spacing of $\left(Y,T^{k}\right)$ for $\delta$.
By appending at most $k-1$ orbit coordinates to the left and to the
right we can obtain two different orbit segments of length $l$ satisfying
$l'\leq l\leq l'+2k-2$: $y^{\left(1\right)},Ty^{\left(1\right)},\dots,T^{\left(l-1\right)k}y^{\left(1\right)}$
and $y^{\left(2\right)},Ty^{\left(2\right)},\dots,T^{\left(l-1\right)k}y^{\left(2\right)}$
- where $y^{\left(1\right)},y^{\left(2\right)}$ belong to the eventual
image of $\left(Y,T^{k}\right)$ - that have an $\varepsilon_{0}$-overlap
and $\varepsilon_{0}$-self-overlaps all smaller than $\frac{l}{3}$
and $l\geq3u_{\delta}$.\\

We claim that these last two orbit segments fulfill for the system
$\left(\stackrel[i=0]{k-1}{\cup}T^{i}Y,T\right)$ what is required
in the remark following the proof of Theorem 4.2.1. The specification
property of $\left(Y,T^{k}\right)$ lets us combine them as we wish
with spaces $ku_{\delta}$ and $\frac{\varepsilon_{0}}{3}$-shadowings
(because of the choice of $\delta$). In addition, they are of length
$lk\geq3ku_{\delta}$. $\blacksquare$\\

Einstein Institute of Mathematics, Edmond J. Safra campus, The Hebrew
University of Jerusalem, Israel.

matan.tal@mail.huji.ac.il
\end{document}